\documentclass[preprint]{elsarticle}
\usepackage{amsmath,amssymb,graphicx}
\newtheorem{df}{Definition}[section]

\newtheorem{thm}[df]{Theorem}
\newtheorem{lem}[df]{Lemma}

\title{On the uncertainty product of spherical functions}
\author[IIN]{I. Iglewska--Nowak\corref{cor}}
\address[IIN]{West Pomeranian University of Technology, Szczecin, School of Mathematics, al.~Piast\'ow~17, 70--310 Szczecin, Poland}
\cortext[cor]{Corresponding author, e--mail: iiglewskanowak@zut.edu.pl, phone/fax: 0048 91 449 4826}

\begin{document}

\begin{abstract}The uncertainty product of a function is a quantity that measures the trade-off between the space and the frequency localization of the function. Its boundedness from below is the content of various uncertainty principles. In the present paper, functions over the $n$-dimensional sphere are considered. A formula is derived that expresses the uncertainty product of a continuous function in terms of its Fourier coefficients. It is applied to a directional derivative of a zonal wavelet, and the behavior of the uncertainty product of this function is discussed.
\end{abstract}

\maketitle

\section{Introduction}

Similarly as in physics (the Heisenberg uncertainty principle), several uncertainty principles are valid in mathematics --- for functions on distinct domains. They state that a function cannot be sharp simultaneously in the space and the momentum domain. This is expressed quantitatively by boundedness from below of an uncertainty constant, defined depending on the function domain.

 The present paper deals with functions over the $n$-dimensional sphere. An uncertainty principle for spherical functions was for the first time formulated and proved in~\cite{NW96} for signals (functions) over the two-dimensional sphere. An extension to higher dimensions appeared in~\cite{RV97}, though only for zonal (rotation-invariant) functions. The result from~\cite{GG04a} is a generalization of the theorems from \cite{NW96,RV97} and it is now regarded as the canonical version of an uncertainty principle on the sphere.

\begin{thm}\label{thm:GG} \emph{\cite[Theorem~4.3]{GG04a}} Let $F$ be a non-zero $\mathcal C^1$-function over the $n$-dimensional sphere~$\mathcal S^n$. Then,
\begin{align*}
&\left[\|F\|_2^2-\frac{\left|\left<\circ\,F(\circ),F(\circ)\right>\right|^2}{\|F\|_2^2}\right]^{1/2}
   \cdot\left[\|\nabla_{\mathcal S^n}F\|_2^2-\frac{\left|\left<\nabla_{\mathcal S^n}F,F\right>\right|^2}{\|F\|_2^2}\right]^{1/2}\\
&\qquad\geq\frac{n}{2}\cdot\frac{\left|\left<\circ\,F(\circ),F(\circ)\right>\right|^2}{\|F\|_2^2}.
\end{align*}
\end{thm}

Further papers yield different proofs for this result resp. its weaker version
\begin{equation}\label{eq:UP}
\left[\|F\|_2^2-\frac{\left|\left<\circ\,F(\circ),F(\circ)\right>\right|^2}{\|F\|_2^2}\right]^{1/2}\cdot\|\nabla_{\mathcal S^n}F\|_2
   \geq\frac{n}{2}\cdot\left|\left<\circ\,F(\circ),F(\circ)\right>\right|,
\end{equation}
see~\cite{GG04b} resp. \cite{DX15a} together with \cite{DX15b}. The term $\|\nabla_{\mathcal S^n}F\|_2$ in equation~\eqref{eq:UP} is equal to and it is often replaced by $\sqrt{\left<-\Delta_{\mathcal S^n}F,F\right>}$. A slightly modified version of Theorem~\ref{thm:GG} is a special case of a general result for compact manifolds~\cite{sS15}.

A different uncertainty principle for functions over the two-dimensional sphere can be obtained with operator approach~\cite{kS02}. A study of an uncertainty principle for signal on the unit sphere in the Clifford algebra setting is the content of the paper~\cite{DQC17}.

The left-hand-side of~\eqref{eq:UP} divided by the module of $\left<\circ\,F(\circ),F(\circ)\right>$ is called the uncertainty product (or constant) of the function~$F$. In a series of publications \cite{LFP02,nLF03,nLF07} Noem\'i La\'in Fern\'andez computes the uncertainty products for certain families of functions and shows their optimality. For this purpose, she derives a series representation of the uncertainty product of a function in terms of its Fourier coefficients, in her doctoral thesis~\cite{nLF03} for zonal functions over the two-dimensional sphere (the result was later generalized by the author of the present paper to the case of $n$-dimensional spheres~\cite{IIN16MR}), and in~\cite{nLF07} for functions over the $n$-dimensional sphere with the gravity center $\int_{\mathcal S^n}x\,|F(x)|^2\,d\sigma(x)$ directed along the $(n+1)^{\text{st}}$ axis. The latter result is sufficient for the purpose of finding of optimal functions with respect to the uncertainty product. However, if a function is given, an orthogonal transformation is required in order to shift the gravity center to the required direction. This orthogonal transform cannot be found unless the actual direction of the gravity center is known.

The purpose of the present paper is to derive a series representation for the center of gravity of a fuction over~$\mathcal S^n$, as well as for the so-called space and the momentum the variances. Further, the formulas are applied to the second directional derivative of the Poisson wavelet~$g_\rho^1$.  This function is an example of directional wavelets, introduced by Hayn and Holschneider in~\cite{HH09} and intended for analyzing of spherical signals with directional features. The result is compared to that obtained for the zonal Poisson wavelets~\cite{IIN17UPW} and discussed in view of the general result for zonal wavelets~\cite{IIN17USW}.

The paper is organized as follows.  Section~\ref{sec:preliminaries} contains basic information about functions over $n$-dimensional spheres. In Section~\ref{sec:UP} series representation of the gravity center, as well as the variances in the space and the momentum domains are derived. Finally, in Section~\ref{sec:DW} directional wavelets are briefly characterized and the formulas obtained previously are applied to a representative of this function family.

\section{Preliminaries}\label{sec:preliminaries}

\subsection{Functions on the sphere}

By $\mathcal{S}^n$ we denote the $n$--dimensional unit sphere in $n+1$--dimensional Euclidean space~$\mathbb{R}^{n+1}$ with the rotation--invariant measure~$d\sigma$ normalized such that
$$
\int_{\mathcal{S}^n}d\sigma=1.
$$
The surface element $d\sigma$ is explicitly given by
$$
d\sigma=\frac{1}{\Sigma_n}\sin^{n-1}\vartheta_1\,\sin^{n-2}\vartheta_2\dots\sin\vartheta_{n-1}d\vartheta_1\,d\vartheta_2\dots d\vartheta_{n-1}d\varphi,
$$
for
$$
\Sigma_n=\frac{2\pi^{\frac{n+1}{2}}}{\Gamma\left(\frac{n+1}{2}\right)},
$$
where $(r,\vartheta_1,\vartheta_2,\dots,\vartheta_{n-1},\varphi)\in[0,\infty)\times[0,\pi]^{n-1}\times[0,2\pi)$ are the spherical coordinates satisfying
\begin{equation*}\begin{split}
x_1&=r\cos\vartheta_1,\\
x_2&=r\sin\vartheta_1\cos\vartheta_2,\\
x_3&=r\sin\vartheta_1\sin\vartheta_2\cos\vartheta_3,\\
&\dots\\
x_{n-1}&=r\sin\vartheta_1\sin\vartheta_2\dots\sin\vartheta_{n-2}\cos\vartheta_{n-1},\\
x_n&=r\sin\vartheta_1\sin\vartheta_2\dots\sin\vartheta_{n-2}\sin\vartheta_{n-1}\cos\varphi,\\
x_{n+1}&=r\sin\vartheta_1\sin\vartheta_2\dots\sin\vartheta_{n-2}\sin\vartheta_{n-1}\sin\varphi.
\end{split}\end{equation*}

The $\mathcal{L}^p(\mathcal{S}^n)$--norm of a function is given by
$$
\|F\|_{\mathcal{L}^p(\mathcal{S}^n)}=\left[\int_{\mathcal S^n}|F(x)|^p\,d\sigma(x)\right]^{1/p}.
$$
The scalar product of $F,G\in\mathcal L^2(\mathcal S^n)$ is defined by
$$
\left<F,G\right>_{\mathcal L^2(\mathcal S^n)}=\int_{\mathcal S^n}\overline{F(x)}\,G(x)\,d\sigma(x),
$$
such that $\|F\|_2^2=\left<F,F\right>$.

Gegenbauer polynomials $C_l^\lambda$ of order~$\lambda\in\mathbb R$ and degree $l\in\mathbb{N}_0$, are defined in terms of their generating function
\begin{equation*}
\sum_{l=0}^\infty C_l^\lambda(t)\,r^l=\frac{1}{(1-2tr+r^2)^\lambda},\qquad t\in[-1,1].
\end{equation*}

A set of Gegenbauer polynomials $\bigl\{C_l^\lambda\bigr\}_{l\in\mathbb N_0}$ builds a complete orthogonal system on $[-1,1]$ with weight $(1-t^2)^{\lambda-1/2}$. Consequently, it is an orthogonal basis for zonal functions on the $(2\lambda+1)$--dimensional sphere. Thus, the number
$$
\lambda=\frac{n-1}{2}
$$
is used interchangeable with~$n$. The following relations are valid for Gegenbauer polynomials:

\begin{align}
2(l+\lambda)\,t\,C_l^\lambda&=(l+2\lambda-1)\,C_{l-1}^\lambda(t)+(l+1)\,C_{l+1}^\lambda(t)\label{eq:rec1}\\
(l+1)\,C_{l+1}^\lambda(t)&=2\lambda\left[t\,C_l^{\lambda+1}(t)-C_{l-1}^{\lambda+1}(t)\right]\label{eq:rec2}\\
\int_{-1}^1C_{l_1}^\lambda(t)\,C_{l_2}^\lambda(t)\,(1-t^2)^{\lambda-\frac{1}{2}}\,dt
   &=\delta_{l_1l_2}\cdot\frac{\pi\,\Gamma(l+2\lambda_1)}{2^{2\lambda-1}\,l_1!\,(l+\lambda_1)\,[\Gamma(\lambda)]^2}\label{eq:SP_Gegenbauer}
\end{align}
cf. \cite{GR}, formulas~8.933.1-2 and~8.939.8.

Let $Q_l$ denote a polynomial on~$\mathbb{R}^{n+1}$ homogeneous of degree~$l$, i.e., such that $Q_l(az)=a^lQ_l(z)$ for all $a\in\mathbb R$ and $z\in\mathbb R^{n+1}$, and harmonic in~$\mathbb{R}^{n+1}$, i.e., satisfying $\nabla^2Q_l(z)=0$, then $Y_l(x)=Q_l(x)$, $x\in\mathcal S^n$, is called a hyperspherical harmonic of degree~$l$. The set of hyperspherical harmonics of degree~$l$ restricted to~$\mathcal S^n$ is denoted by $\mathcal H_l(\mathcal S^n)$. Hyperspherical harmonics of distinct degrees are orthogonal to each other. The number of linearly independent hyperspherical harmonics of degree~$l$ is equal to
$$
N=\frac{2(l+\lambda)(l+2\lambda-1)!}{l!\,(2\lambda)!}.
$$

The hyperspherical harmonics are eigenfunctions of the operator~$\Delta^\ast$,
\begin{equation}\label{eq:HSH_eigenfunctions}
\Delta^\ast Y(x)=-l\,(l+2\lambda)\,Y(x),\qquad\text{for }Y\in\mathcal H_l(\mathcal S^n),
\end{equation}
where the spherical Laplace-Beltrami operator~$\Delta^\ast$ is the tangential part of the Laplace operator,
$$
\Delta=\frac{\partial^2}{\partial r^2}+\frac{n}{r}\frac{\partial}{\partial r}+\frac{1}{r^2}\Delta^\ast.
$$

In this paper, we will be working with the orthonormal basis for~$\mathcal L^2(\mathcal S^n)=\overline{\bigoplus_{l=0}^\infty\mathcal H_l}$, consisting of hyperspherical harmonics given by
\begin{equation*}
Y_l^k(x)=A_l^k\prod_{\tau=1}^{n-1}C_{k_{\tau-1}-k_\tau}^{\frac{n-\tau}{2}+k_\tau}(\cos\vartheta_\tau)\sin^{k_\tau}\!\vartheta_\tau\cdot e^{k_{n-1}\varphi}
\end{equation*}
with $l=k_0\geq k_1\geq\dots\geq|k_{n-1}|\geq0$, $k$ being a sequence $(k_1,\dots,k_{n-1})$ of integer numbers, and normalization constants
\begin{equation}\label{eq:Alk}
A_l^k=\left(\frac{1}{\Gamma\left(\frac{n+1}{2}\right)}\prod_{\tau=1}^{n-1}\frac{2^{n-\tau+2k_\tau-2}\,(k_{\tau-1}-k_\tau)!\,(n-\tau+2k_{\tau-1})\,\Gamma^2(\frac{n-\tau}{2}+k_\tau)}{\sqrt\pi\,(n-\tau+k_{\tau-1}+k_\tau-1)!}\right)^{1/2},
\end{equation}
compare~\cite[Sec. IX.3.6, formulas (4) and (5)]{Vilenkin}. The set of non-increasing sequences~$k$ in $\mathbb N_0^{n-1}\times\mathbb Z$ with elements bounded by~$l$ will be denoted by $\mathcal M_n(l)$.

Any function $F\in\mathcal L^2(\mathcal S^n)$ has a unique representation as a mean--convergent series
\begin{equation*}
F(x)=\sum_{l=0}^\infty\sum_{k\in\mathcal M_{n-1}(l)} \widehat F_l^k\,Y_l^k(x),\qquad x\in\mathcal S^n,
\end{equation*}
where
$$
\widehat F_l^k=\int_{\mathcal S^n}\overline{Y_l^k(x)}\,F(x)\,d\sigma(x)=\left<Y_l^k,F\right>,
$$
for proof cf.~\cite{Vilenkin}. $\widehat F_l^k$, $l\in\mathbb N_0$, $k\in\mathcal M_n(l)$, are the Fourier coefficients of the function~$F$.

\section{The uncertainty product}\label{sec:UP}

\begin{df}
Suppose, $F\in\mathcal C(\mathcal S^2)$ is a nonzero function. The quantity
\begin{align}\label{eq:xi_O_F}
\xi_O(F)&=\frac{1}{\|F\|_2^2}\int_{\mathcal S^n}x\,|F(x)|^2\,d\sigma(x)
\end{align}
is called its center of gravity in the space domain.
\end{df}

\begin{df}
The the variances in the space and the momentum domain of a $\mathcal C^2(\mathcal S^n)$--function~$F$ with $\xi_O(F)\ne0$ are given by
\begin{equation*}
\text{var}_S(F)=\frac{1-\|\xi_O(F)\|^2}{\|\xi_O(F)\|^2}
\end{equation*}
and
\begin{equation}\label{eq:var_M_F}
\text{var}_M(F)=-\frac{1}{\|F\|_2^2}\int_{\mathcal S^n}\Delta^\ast F(x)\cdot \bar F(x)\,d\sigma(x),
\end{equation}
respectively. The quantity
$$
U(F)=\sqrt{\text{var}_S(F)}\cdot\sqrt{\text{var}_M(F)}
$$
is called the uncertainty product (or constant) of~$F$.
\end{df}

With this notation, the uncertainty principle~\eqref{eq:UP} can be expressed as
$$
U(F)\geq\frac{n}{2}.
$$

In the following, let~$F$ be given by its Fourier expansion,
\begin{equation}\label{eq:F_Fourier_series}
F=\sum_{l=0}^\infty\sum_{k\in\mathcal M_{n-1}(l)}\widehat F_l^k\cdot Y_l^k.
\end{equation}

\subsection{Localization in the space domain}

Substitute~\eqref{eq:F_Fourier_series} into~\eqref{eq:xi_O_F} to see that the center of gravity in the space domain of the function~$F$ is given by
\begin{align}
\|F\|_2^2\cdot\xi_O(F)&=\int_{\mathcal S^n}x\cdot\sum_{l_1,k}\overline{\widehat F_{l_1}^k\,Y_{l_1}^k(x)}
   \cdot\sum_{l_2,m}\widehat F_{l_2}^m\,Y_{l_2}^m(x)\,d\sigma(x)\notag\\
&=\sum_{l_1,k,l_2,m}\overline{\widehat F_{l_1}^k}\cdot\widehat F_{l_2}^m\cdot I(l_1,k,l_2,m)\label{eq:xi_O_F_series}
\end{align}
for
$$
I(l_1,k,l_2,m):=\int_{\mathcal S^n}x\cdot\overline{Y_{l_1}^k(x)}\cdot Y_{l_2}^m(x)\,d\sigma(x).
$$
The notation will be simplified to $I(\widehat k, \widehat m)$ with $\widehat k=(k_0,k_1,\dots,k_{n-1})$, $k_0=l_1$, and $\widehat m=(m_0,m_1,\dots,m_{n-1})$, $m_0=l_2$.

\subsubsection*{Computation of $I(\widehat k,\widehat m)$}

\begin{equation}\label{eq:I}\begin{split}
I&=I(\widehat k,\widehat m)\\
&=\frac{A_{k_0}^k\cdot A_{m_0}^m}{\Sigma_n}\cdot\int_0^\pi\int_0^\pi\cdots\int_0^\pi\int_0^\pi\int_0^{2\pi}
   \prod_{\tau=1}^{n-1}C_{k_{\tau-1}-k_\tau}^{\frac{n-\tau}{2}+k_\tau}(\cos\vartheta_\tau)\cdot\sin^{k_\tau}\vartheta_\tau\cdot e^{-ik_{n-1}\varphi}\\
&\qquad\cdot\prod_{\iota=1}^{n-1}C_{m_{\iota-1}-m_\iota}^{\frac{n-\iota}{2}+m_\iota}(\cos\vartheta_\iota)\cdot\sin^{m_\iota}\vartheta_\iota\cdot e^{im_{n-1}\varphi}\\
&\qquad\cdot\left(\begin{array}{l}
   \cos\vartheta_1\\
   \sin\vartheta_1\,\cos\vartheta_2\\
   \sin\vartheta_1\,\sin\vartheta_2\,\cos\vartheta_3\\
   \vdots\\
   \sin\vartheta_1\,\sin\vartheta_2\cdots\cdots\sin\vartheta_{n-2}\,\cos\vartheta_{n-1}\\
   \sin\vartheta_1\,\sin\vartheta_2\cdots\cdots\sin\vartheta_{n-2}\,\sin\vartheta_{n-1}\,\cos\varphi\\
   \sin\vartheta_1\,\sin\vartheta_2\cdots\cdots\sin\vartheta_{n-2}\,\sin\vartheta_{n-1}\,\sin\varphi
\end{array}\right)d\varphi\\
&\qquad\cdot\sin\vartheta_{n-1}\,d\vartheta_{n-1}\cdot\sin^2\vartheta_{n-2}\,d\vartheta_{n-2}\cdots\sin^{n-2}\vartheta_2\,d\vartheta_2\cdot\sin^{n-1}\vartheta_1\,d\vartheta_1.
\end{split}\end{equation}

The components of~$I$ are products of the following single integrals:
\begin{align*}
C_{\vartheta,1}&(\iota,k_{\iota-1},k_\iota,m_{\iota-1},m_\iota)\\
&:=\int_0^\pi C_{k_{\iota-1}-k_\iota}^{\frac{n-\iota}{2}+k_\iota}(\cos\vartheta_\iota)
   \cdot C_{m_{\iota-1}-m_\iota}^{\frac{n-\iota}{2}+m_\iota}(\cos\vartheta_\iota)\cdot\sin^{k_\iota+m_\iota+n-\iota}\,d\vartheta_\iota,\\
C_{\vartheta,c}&(\iota,k_{\iota-1},k_\iota,m_{\iota-1},m_\iota)\\
   &:=\int_0^\pi C_{k_{\iota-1}-k_\iota}^{\frac{n-\iota}{2}+k_\iota}(\cos\vartheta_\iota)\cdot C_{m_{\iota-1}-m_\iota}^{\frac{n-\iota}{2}+m_\iota}(\cos\vartheta_\iota)
   \cdot\cos\vartheta_\iota\cdot\sin^{k_\iota+m_\iota+n-\iota}\,d\vartheta_\iota,\\
C_{\vartheta,s}&(\iota,k_{\iota-1},k_\iota,m_{\iota-1},m_\iota)\\
   &:=\int_0^\pi C_{k_{\iota-1}-k_\iota}^{\frac{n-\iota}{2}+k_\iota}(\cos\vartheta_\iota)
   \cdot C_{m_{\iota-1}-m_\iota}^{\frac{n-\iota}{2}+m_\iota}(\cos\vartheta_\iota)\cdot\sin^{k_\iota+m_\iota+n-\iota+1}\,d\vartheta_\iota,\\
C_{\varphi,1}&(m_{n-1}-k_{n-1}):=\int_0^{2\pi}e^{i(m_{n-1}-k_{n-1})\varphi}\,d\varphi,\\
C_{\varphi,c}&(m_{n-1}-k_{n-1}):=\int_0^{2\pi}e^{i(m_{n-1}-k_{n-1})\varphi}\,\cos\varphi\,d\varphi,\\
C_{\varphi,s}&(m_{n-1}-k_{n-1}):=\int_0^{2\pi}e^{i(m_{n-1}-k_{n-1})\varphi}\,\sin\varphi\,d\varphi.
\end{align*}

More exactly,
\begin{equation}\label{eq:I_tau}
I_\tau=\frac{A_{l_1}^k\cdot A_{l_2}^m}{\Sigma_n}
   \cdot\prod_{\iota=1}^{\tau-1}C_{\vartheta,s}^\iota\cdot C_{\vartheta,c}^\tau\cdot\prod_{\iota=\tau+1}^{n-1}C_{\vartheta,1}^\iota\cdot C_{\varphi,1}.
\end{equation}
for $\tau=1,2,\dots,n-1$ and
\begin{align}
I_n&=\frac{A_{l_1}^k\cdot A_{l_2}^m}{\Sigma_n}\cdot\prod_{\iota=1}^{n-1}C_{\vartheta,s}^\iota\cdot C_{\varphi,c},\label{eq:I_n}\\
I_{n+1}&=\frac{A_{l_1}^k\cdot A_{l_2}^m}{\Sigma_n}\cdot\prod_{\iota=1}^{n-1}C_{\vartheta,s}^\iota\cdot C_{\varphi,s},\label{eq:I_{n+1}}
\end{align}
where the superscript~$\iota$ is used instead of the set of variables $(\iota,k_{\iota-1},k_\iota,m_{\iota-1},m_\iota)$ and the variables of $I_\tau$  and $C_\varphi$ are omitted.\\

For $C_\vartheta$ use~\eqref{eq:SP_Gegenbauer} as well as the following two lemmas.

\begin{lem}\label{lem:C_theta_c}
The integral
$$
\int_0^\pi t\cdot C_{l_1}^\lambda(t)\cdot C_{l_2}^\lambda(t)\cdot(1-t^2)^{\lambda-\frac{1}{2}}\,dt
$$
disappears for $l_2-l_1\ne\pm1$. Otherwise,
\begin{equation*}
\int_0^\pi t\cdot C_{l}^\lambda(t)\cdot C_{l+1}^\lambda(t)\cdot(1-t^2)^{\lambda-\frac{1}{2}}\,dt
   =\frac{\pi\,\Gamma(l+2\lambda+1)}{2^{2\lambda}\,l!\,(l+\lambda)\,(l+\lambda+1)\,[\Gamma(\lambda)]^2}.
\end{equation*}
\end{lem}

{\bf Proof. }According to~\eqref{eq:rec1}, $t\cdot C_{l_1}^\lambda(t)$ can be expressed as
$$
t\,C_{l_1}^\lambda(t)=\frac{l+2\lambda_1-1}{2(l+\lambda_1)}\,C_{l_1-1}^\lambda(t)+\frac{l_1+1}{2(l+\lambda_1)}\,C_{l_1+1}^\lambda(t).
$$
Thus, by~\eqref{eq:SP_Gegenbauer},
\begin{align*}
\int_0^\pi &t\cdot C_{l_1}^\lambda(t)\cdot C_{l_2}^\lambda(t)\cdot(1-t^2)^{\lambda-\frac{1}{2}}\,dt\\
&=\frac{\pi\,\Gamma(l+2\lambda_1)}{2^{2\lambda}\,(l_1-1)!\,(l+\lambda_1-1)\,(l+\lambda_1)\,[\Gamma(\lambda)]^2}\cdot\delta_{l_1-1,l_2}\\
&+\frac{\pi\,\Gamma(l+2\lambda_1+1)}{2^{2\lambda}\,l_1!\,(l+\lambda_1)\,(l+\lambda_1+1)\,[\Gamma(\lambda)]^2}\cdot\delta_{l_1+1,l_2}\\
&=\frac{\pi\,\Gamma(l+2\lambda_1)}{2^{2\lambda}\,l_2!\,(l+\lambda_2)\,(l+\lambda_1)\,[\Gamma(\lambda)]^2}\cdot\delta_{l_1-1,l_2}\\
&+\frac{\pi\,\Gamma(l+2\lambda_2)}{2^{2\lambda}\,l_1!\,(l+\lambda_1)\,(l+\lambda_2)\,[\Gamma(\lambda)]^2}\cdot\delta_{l_1+1,l_2}.
\end{align*}
\hfill$\Box$

\begin{lem}\label{lem:C_theta_s}
\begin{align*}
\int_0^\pi &C_{l_1}^{\lambda}(t)\cdot C_{l_2}^{\lambda+1}(t)\cdot(1-t^2)^{\lambda+\frac{1}{2}}\,dt\\
&=\frac{\pi\,\lambda\,\Gamma(l+2\lambda_1+2)}{2^{2\lambda+1}\,l_1!\,(l+\lambda_1)\,(l+\lambda_1+1)\,[\Gamma(\lambda+1)]^2}\cdot\delta_{l_1l_2}\\
&-\frac{\pi\,\lambda\,\Gamma(l+2\lambda_1)}{2^{2\lambda+1}\,(l_1-2)!\,(l+\lambda_1-1)\,(l+\lambda_1)\,[\Gamma(\lambda+1)]^2}\cdot\delta_{l_1-2,l_2}.
\end{align*}
\end{lem}

{\bf Proof. }According to~\eqref{eq:rec2}, $C_{l_1}^{\lambda}$ can be written as
$$
C_{l_1}^{\lambda}(t)=\frac{2\lambda}{l_1}\left[t\,C_{l_1-1}^{\lambda+1}(t)-C_{l_1-2}^{\lambda+1}(t)\right]
$$
and further, by~\eqref{eq:rec1},
\begin{align*}
C_{l_1}^{\lambda}(t)&=\frac{2\lambda}{l_1}\left[\frac{(l+2\lambda_1)\,C_{l_1-2}^{\lambda+1}(t)+l_1C_{l_1}^{\lambda+1}(t)}{2(l+\lambda_1)}-C_{l_1-2}^{\lambda+1}(t)\right]\\
&=\frac{\lambda}{l+\lambda_1}\left[C_{l_1}^{\lambda+1}(t)-C_{l_1-2}^{\lambda+1}(t)\right].
\end{align*}
The assertion follows  by~\eqref{eq:SP_Gegenbauer}.\hfill$\Box$\\

Summarizing,
\begin{equation*}
C_{\vartheta,1}(\iota,a,b,a,b)=\frac{\pi\,\Gamma(n-\iota+a+b)}{2^{n-\iota+2b-1}\,(a-b)!\,(\frac{n-\iota}{2}+a)\,[\Gamma(\frac{n-\iota}{2}+b)]^2}
\end{equation*}
and 
\begin{equation*}
C_{\vartheta,1}(\iota,a,b,c,d)=0
\end{equation*}
if $a\ne c$ or $b\ne d$. Further,
\begin{align*}
C_{\vartheta,c}&(\iota,a,b,a+1,b)=C_{\vartheta,c}(\iota,a+1,b,a,b)\\
&=\frac{\pi\,\Gamma(n-\iota+a+b+1)}{2^{n-\iota+2b}\,(a-b)!\,(\frac{n-\iota}{2}+a)\,(\frac{n-\iota}{2}+a+1)\,[\Gamma(\frac{n-\iota}{2}+b)]^2}
\end{align*}
and
\begin{equation*}
C_{\vartheta,c}(\iota,a,b,c,b)=0
\end{equation*}
if $a-c\ne\pm1$. Finally,
\begin{align*}
C_{\vartheta,s}&(\iota,a,b,a+1,b+1)=C_{\vartheta,s}(\iota,a+1,b+1,a,b)\\
&=\frac{\pi\,(\frac{n-\iota}{2}+b)\,\Gamma(n-\iota+a+b+2)}{2^{n-\iota+2b+1}\,(a-b)!\,(\frac{n-\iota}{2}+a)\,(\frac{n-\iota}{2}+a+1)\,[\Gamma(\frac{n-\iota}{2}+b+1)]^2}\\
C_{\vartheta,s}&(\iota,a,b+1,a+1,b)=C_{\vartheta,s}(\iota,a+1,b,a,b+1)\\
&=-\frac{\pi\,(\frac{n-\iota}{2}+b)\,\Gamma(n-\iota+a+b+1)}{2^{n-\iota+2b+1}\,(a-b-1)!\,(\frac{n-\iota}{2}+a)\,(\frac{n-\iota}{2}+a+1)\,[\Gamma(\frac{n-\iota}{2}+b+1)]^2}
\end{align*}
and
\begin{equation*}
C_{\vartheta,s}(\iota,a,b,c,d)=0
\end{equation*}
if $b-d\ne\pm1$ or $a-c\ne\pm1$.\\

The values of~$C_\vartheta$ for other sets of arguments are not needed in the computation of the space localization of a function. It is a consequence of Theorem~\ref{thm:space_localization_HSH}.

$C_\varphi$ can be obtained by a simple computation. They are given by
\begin{equation}\label{eq:C_phi_1}
C_{\varphi,1}(a)=\begin{cases}2\pi&\text{for }a=0,\\0&\text{otherwise,}\end{cases}
\end{equation}
and
\begin{equation}\label{eq:C_phi_c}
C_{\varphi,c}(a)=\begin{cases}\pi&\text{for }a=\pm1,\\0&\text{otherwise,}\end{cases}
\end{equation}
as well as
\begin{equation}\label{eq:C_phi_s}
C_{\varphi,s}(a)=\begin{cases}ia\pi&\text{for }a=\pm1,\\0&\text{otherwise.}\end{cases}
\end{equation}

After this preparation we are able to prove the following theorem.

\begin{thm}\label{thm:space_localization_HSH}If $I(\widehat k,\widehat m)\ne0$ then there exists $\nu\in\{1,2,\dots,n\}$ such that $\widehat k$ and $\widehat m$ satisfy
\begin{equation}\label{eq:k_m_structure}\begin{split}
m_\iota=k_\iota\pm1\qquad&\text{for }\iota<\nu,\\
m_\iota=k_\iota\hspace{3.7em}&\text{for }\iota\geq\nu.
\end{split}\end{equation}
If $\nu<n$, then $I_\nu\ne0$ and $I_\iota=0$ for $\iota\ne\nu$. If $\nu=n$, then $I_n\ne0$, $I_{n+1}\ne0$, and $I_\iota=0$ for $\iota<\nu$. Non-vanishing components of~$I$ are given by~\eqref{eq:I_tau} or~\eqref{eq:I_n} together with~\eqref{eq:I_{n+1}}.
\end{thm}

{\bf Proof. }The structure of the components of~$I$, given by~\eqref{eq:I_tau}, \eqref{eq:I_n} and~\eqref{eq:I_{n+1}} is obtained directly from~\eqref{eq:I}.\\

{\bf Case 1:} Suppose $I_\tau\ne0$ for $\tau<n$.\\ 

Then, by~\eqref{eq:I_tau}, $C_{\varphi,1}\ne0$. It follows from~\eqref{eq:C_phi_1} that $m_{n-1}=k_{n-1}$.\\

Further, $C_{\vartheta,1}^\iota\ne0$ for $\iota=\tau+1,\tau+2,\dots,n-1$. Using this fact, we shall show by induction that $m_\iota=k_\iota$ for $\iota=n-2,n-3,\dots,\tau$. For $\iota=n-1,n-2,\dots,\tau+1$ suppose, $m_\iota=k_\iota$ (the assumption is satisfied for $\iota=n-1$ by the previous consideration). Then, in the integral representing~$C_{\vartheta,1}^\iota$, the Gegenbauer polynomials have the same order. According to~\eqref{eq:SP_Gegenbauer}, the degrees of the polynomials are equal to each other. Thus, $m_{\iota-1}=k_{\iota-1}$.\\

$C_{\vartheta,c}^\tau\ne0$ and $m_\tau=k_\tau$ imply $m_{\tau-1}-k_{\tau-1}=\pm1$, see~\eqref{eq:C_theta_c} and Lemma~\ref{lem:C_theta_c}.\\

Finally, $C_{\vartheta,s}^\iota\ne0$ for $\iota=1,2,\dots,\tau-1$. Thus, by induction, $m_\iota-k_\iota=\pm1$ for $\iota=\tau-2,\tau-3,\dots,1$. More exactly, suppose for $\iota=\tau-1,\tau-2,\dots,2$ that $m_\iota-k_\iota=\pm1$ (this condition is satisfied for $\iota=\tau-1$). In the case $m_\iota-k_\iota=1$ it follows from~\eqref{eq:C_theta_s} and Lemma~\ref{lem:C_theta_s} that
$$
m_{\iota-1}-m_\iota=k_{\iota-1}-k_\iota\qquad\lor\qquad m_{\iota-1}-m_\iota=k_{\iota-1}-k_\iota-2.
$$
Consequently,
$$
m_{\iota-1}-k_{\iota-1}=\pm1.
$$
For $m_\iota-k_\iota=-1$ change the roles of $k$ and $m$.\\

Summarizing, \eqref{eq:k_m_structure} is satisfied with $\nu=\tau$. Thus, $I_\iota=0$ for $\iota\ne\tau$, $\iota<n$. (Otherwise, if $I_\iota\ne0$, then~\eqref{eq:k_m_structure} would be satisfied for $\nu=\iota\ne\tau$.) Further, $I_n=I_{n+1}=0$ since $C_{\varphi,c}=C_{\varphi,s}=0$ for $m_{n-1}=k_{n-1}$.\\

{\bf Case 2:} Suppose, $I_n\ne0$ or $I_{n+1}\ne0$.\\

Then, $C_{\varphi,c}\ne0$. Consequently, by~\eqref{eq:C_phi_c}, $m_{n-1}-k_{n-1}=\pm1$ and therefore \mbox{$C_{\varphi,1}=0$}. This implies $I_\iota=0$ for $\iota<n$. It follows from $C_{\varphi,s}\ne0$ and \mbox{$m_{n-1}-k_{n-1}=\pm1$} by induction that $m_\iota-k_\iota=\pm1$ for \mbox{$\iota=n-2,n-3,\dots,1$}.\hfill$\Box$\\

\subsection*{Computation of~$\xi_O(F)$ and var$_S(F)$}

\begin{df}Let $\widehat k=(k_0,k_1,\dots,k_{n-1})$ and $\widehat m=(m_0,m_1,\dots,m_{n-1})$ be sequences in $\mathbb N_0^n\times\mathbb Z$ such that $(k_1,k_2,\dots,k_{n-1})\in\mathcal M_{n-1}(k_0)$ and $(m_1,m_2,\dots,m_{n-1})\in\mathcal M_{n-1}(m_0)$. We call $\widehat k$ and $\widehat m$ $\nu$-conjugate if they satisfy~\eqref{eq:k_m_structure} for a constant $\nu\in\{1,2,\dots,n\}$.
\end{df}

In view of Theorem~\ref{thm:space_localization_HSH}, the summation in~\eqref{eq:xi_O_F_series} goes over $\nu$-conjugate sequences~$\widehat k$ and~$\widehat m$ for $\nu=1,2,\dots,n$.

\begin{thm} Let~$F$ be a continuous functions over~$\mathcal S^n$ with $\|F\|\ne0$. Then,
$$
\|F\|_2^2\cdot\xi_O(F)
   =\sum_{\nu=1}^n\sum_{\genfrac{}{}{0pt}{1}{\widehat k,\widehat m}{\nu\text{-conjugate}}}\overline{\widehat F_{k_0}^k}\cdot\widehat F_{m_0}^m\cdot I(\widehat k,\widehat m).
$$
Components with distinct $\nu$ are orthogonal to each other. The values of the non-vanishing component(s) of~$I$ are given by~\eqref{eq:I_tau} for $\nu<n$ and \eqref{eq:I_n} together with~\eqref{eq:I_{n+1}} for $\nu=n$. 
\end{thm}

By computation of the space the variance of~$F$ note that by~\eqref{eq:C_phi_c} and~\eqref{eq:C_phi_s},
$$
\left|C_{\varphi,c}\right|^2+\left|C_{\varphi,s}\right|^2=2\pi^2.
$$
Hence, $|I_n|^2+|I_{n+1}|^2=2|I_n|^2$. Further, for each~$\nu$ symmetries in the values of~$I$ with respect to the sign of $m_\iota-k_\iota$ can be exploited in order to simplify the calculation, see, e.g., Subsection~\ref{subs:localization}.

\subsection{Variance in the momentum domain}

The the variance in the momentum domain can be calculated directly from the definition~\eqref{eq:var_M_F} using the eigenfunction property of the hyperspherical harmonics~\eqref{eq:HSH_eigenfunctions} and the orthonormality of~$Y_l^k$.

\begin{thm} For $F\in\mathcal C(\mathcal S^n)$,
\begin{equation}\label{eq:varMF_series}
\text{var}_M(F)=\sum_{l=0}^\infty l(l+2\lambda)\sum_{k\in\mathcal M_{n-1}(l)}\left|\widehat F_l^k\right|^2.
\end{equation}
\end{thm}

\section{Example: directional wavelets}\label{sec:DW}

The idea to construct a not rotation invariant spherical wavelet by a directional derivative of the Poisson kernel was described in~\cite{HH09}. In this case, 'directional derivative' means that the source of the field, or the defining position, is rotated around an axis that is perpendicular to the symmetry axis, and the derivative of this transform is computed. Functions constructed in that way are wavelets according to the definition given in~\cite{mH96}. This concept -- generalized to the $n$-dimensional case -- was also studied by the author of the present paper. I  showed in~\cite{IIN18DW} that directional derivatives of some zonal wavelets satisfy a slightly modified definition of wavelets derived from approximate identities (see \cite{FW96,FW-C97,FGS-book,sB09,EBCK09,BE10WS3,BE10KBW} for the origins of this concept and \cite{IIN15CWT} for a comprehensive survey). In \cite{IIN17FDW} I proposed further relaxation on the constraints on wavelets derived from approximate identities and showed that directional derivatives of a wide class of functions are wavelets according to the new definition.

In \cite[Theorem~4.3]{IIN18DW} a recipe is given how to compute the Fourier coefficients of directional derivatives of a zonal function. We shall apply it to the Poisson wavelet~$g_\rho^1$ over~$\mathcal S^n$, $n=2\lambda+1$ (see~\cite{CPMDHJ05,HI07,IIN15PW} for details on the Poisson wavelets). Consider the second directional derivative of~$g_\rho^1$,
\begin{align*}
G:=\left(g_\rho^1\right)^{[2]}=&\frac{1}{\Sigma_n}\cdot\left[\sum_{l=1}^\infty\beta_{l,0}^2\cdot\frac{l+\lambda}{\lambda}\,\rho l\,e^{-\rho l}
   \cdot\frac{Y_l^{(0,0,\dots,0)}}{A_l^{(0,0,\dots,0)}}\right.\\
&\left.-\sum_{l=2}^\infty\beta_{l,0}\,\beta_{l,1}\cdot\frac{l+\lambda}{\lambda}\,\rho l\,e^{-\rho l}\cdot\frac{Y_l^{(2,0,0,\dots,0)}}{A_l^{(0,0,\dots,0)}}\right],
\end{align*}
where
$$
\beta_{l,k_1}=\sqrt\frac{(k_1+1)(2\lambda+k_1-1)(l-k)(l+2\lambda+k_1)}{(2\lambda+2k_1-1)(2\lambda+2k_1+1)},
$$
see \cite[formula~(20)]{IIN18DW}. (The expression~(55) in~\cite{IIN18DW} contains a mistake: the signs on the right-hand-side should be changed to the opposite ones, as it follows from \cite[Theorem~4.3]{IIN18DW}.)
In order to simplify the notation, we shall abbreviate the indices appearing in the above formula and write $A_l^0$ resp. $Y_l^{k_1}$ instead of $A_l^{(0,0,\dots,0)}$ resp. $Y_l^{(k_1,0,0,\dots,0)}$. Similarly, the non-vanishing Fourier coefficients of~$G$ will be denoted by~$\widehat G_l^0$ and $\widehat G_l^2$. Formula~\eqref{eq:Alk} yields
$$
A_l^0=\sqrt\frac{(2\lambda-1)!\,l!\,(l+\lambda)}{\lambda\,(l+2\lambda-1)}.
$$
The norm of the considered function is given by
\begin{align*}
\left\|G\right\|_2^2&=\frac{1}{\lambda^2\,\Sigma_n^2}
   \sum_{l=1}^\infty\frac{\beta_{l,0}^2\left(\beta_{l,0}^2+\beta_{l,1}^2\right)(l+\lambda)^2}{(A_l^0)^2}\,(\rho l)^2\,e^{-2\rho l},
\end{align*}
i.e.,
\begin{equation}\label{eq:norm_g12_as_series}\begin{split}
\left\|G\right\|_2^2&=\frac{2\rho^2}{\Sigma_n^2\,(2\lambda+1)(2\lambda+3)}\\
&\cdot\sum_{l=1}^\infty\left[3 l^3+9\lambda l^2+2\lambda(3\lambda-2)l-4\lambda^2\right]\binom{l+2\lambda}{l}\,l^3\,e^{-2\rho l}.
\end{split}\end{equation}
Similarly as in~\cite{IIN17UPW} (see the second formula on page~352) we shall introduce the following quantity
\begin{equation}\label{eq:Smmu}
S_m^\mu(\rho)=\sum_{l=0}^\infty\binom{l+\mu}{l}l^m\,e^{-2\rho l}.
\end{equation}
According to \cite[formula~(6)]{IIN17UPW}, for $\mu\geq3$ it behaves like
\begin{equation}\label{eq:Smmurho}\begin{split}
S_m^\mu(\rho)&=\frac{1}{2^{\mu+m+2}}\cdot\left[\frac{2(\mu+m)!}{\mu!\,\rho^{\mu+m+1}}+\frac{2(\mu+1)(\mu+m-1)!}{(\mu-1)!\,\rho^{\mu+m}}\right.\\
&\left.+\frac{(\mu+\tfrac{2}{3})(\mu+1)(\mu+m-2)!}{(\mu-2)!\,\rho^{\mu+m-1}}+\frac{\tfrac{1}{3}\mu(\mu+1)^2(\mu+m-3)!}{(\mu-3)!\,\rho^{\mu+m-2}}\right]\\
&+\mathcal{O}(\rho^{-\mu-m+3})\qquad\text{for }\rho\to0.
\end{split}\end{equation}
For values of~$\mu$ less than~$3$, another asymptotic can be found in~\cite{IIN17UPW}. Here, we shall concentrate only on the one case. The squared norm of~$G$ given by~\eqref{eq:norm_g12_as_series} can be expressed as
\begin{equation}\label{eq:norm2_g12}\begin{split}
\left\|G\right\|_2^2=&\frac{2\rho^2}{\Sigma_n^2\,(2\lambda+1)(2\lambda+3)}\left[3S_6^{2\lambda}+9\lambda S_5^{2\lambda}
   +2\lambda(3\lambda-2)S_4^{2\lambda}-4\lambda^2S_3^{2\lambda}\right].
\end{split}\end{equation}
Substituting~\eqref{eq:Smmurho} into~\eqref{eq:norm2_g12} we obtain
\begin{equation}\label{eq:G2}\begin{split}
\left\|G\right\|_2^2=&\,\frac{\lambda+1}{\Sigma_n^2\rho^{2\lambda}}
   \cdot\left[\frac{3(\lambda+2)(\lambda+3)(2\lambda+5)}{2^{2\lambda+3}\rho^5}+\frac{3\lambda(\lambda+2)^2(2\lambda+5)}{2^{2\lambda+2}\rho^4}\right.\\
&+\frac{\lambda(\lambda+2)(2\lambda+3)(6\lambda^2+11\lambda-3)}{2^{2\lambda+3}\rho^3}\\
&\left.+\frac{\lambda^2(\lambda+1)(2\lambda+3)(2\lambda^2+3\lambda-3)}{2^{2\lambda+2}\rho^2}+\mathcal O\left(\frac{1}{\rho}\right)\right].
\end{split}\end{equation}

\subsection{Localization in the space domain}\label{subs:localization}

In order to find the gravity center of~$G$, note that the only $\nu$-conjugate sequences~$\widehat k$, $\widehat m$ with non-vanishing product $\overline{\widehat{G}_{k_0}^k}\cdot\widehat{G}_{m_0}^m$ are those with
$$
m_0=k_0\pm1,\qquad m=k=(0,0,\dots,0)\lor m=k=(2,0,0,\dots,0).
$$
Consequently, by Theorem~\ref{thm:space_localization_HSH}, the gravity center is directed along the $x_1$-axis.

(Note that by the structure theorem for directional derivatives of zonal functions \cite[Theorem~4.3]{IIN18DW}, the non-vanishing Fourier coefficients~$\widehat F_{k_0}^k$ appearing in the series representation of such a directional derivative~$F$ satisfy either $k_1\in2\mathbb N$ or $k_1\in2\mathbb N+1$ and $k_\iota=0$ for $\iota>1$. Therefore, the only $\nu$-conjugate sequences~$\widehat k$, $\widehat m$ with non-vanishing product $\overline{\widehat{F}_{k_0}^k}\cdot\widehat{F}_{m_0}^m$ are those with $\nu=1$, and the center of gravity of such a function is necessarily directed along the $x_1$-axis.)

Since the Fourier coefficients of~$G$ are real,
\begin{equation}\label{eq:G2_xi_O_G}\begin{split}
\|G\|^2\cdot\xi_O(G)&=2\cdot\sum_{l=1}^\infty\widehat G_l^0\cdot\widehat G_{l+1}^0\cdot I(l,0,l+1,0)\\
&+2\cdot\sum_{l=2}^\infty\widehat G_l^2\cdot\widehat G_{l+1}^2\cdot I(l,2,l+1,2).
\end{split}\end{equation}

In this formula, the second and the fourth arguments of~$I$, denoted by a number~$k_1$, are to be interpreted as sequences $(k_1,0,0,\dots,0)$. Now,
\begin{align*}
\widehat G_l^0&\cdot\widehat G_{l+1}^0
   =\frac{\rho^2\,\beta_{l,0}^2\,\beta_{l+1,0}^2(l+\lambda)\,(l+\lambda+1)\,l\,(l+1)\,e^{-\rho(2l+1)}}{\Sigma_n^2\,\lambda^2\,A_l^0\,A_{l+1}^0}\\
&=\frac{\rho^2}{\Sigma_n^2A_l^0A_{l+1}^0}
   \cdot\frac{l^2(l+1)^2(l+\lambda)(l+\lambda+1)(l+2\lambda)(l+2\lambda+1)}{\lambda^2(2\lambda+1)^2}\cdot e^{-\rho(2l+1)}
\end{align*}
and
\begin{align*}
\widehat G_l^2\cdot\widehat G_{l+1}^2
   =&\frac{\rho^2\,\beta_{l,0}\,\beta_{l,1}\,\beta_{l+1,0}\,\beta_{l+1,1}\,(l+\lambda)\,(l+\lambda+1)\,l\,(l+1)\,e^{-\rho(2l+1)}}{\Sigma_n^2\,\lambda^2\,A_l^0\,A_{l+1}^0}\\
=&\frac{\rho^2}{\Sigma_n^2A_l^0A_{l+1}^0}\cdot\frac{4l^2(l+1)(l+\lambda)(l+\lambda+1)(l+2\lambda+1)}{\lambda(2\lambda+1)(4\lambda^2+8\lambda+3)}\\
   &\cdot\sqrt{(l+2\lambda)(l+2\lambda+2)(l-1)(l+1)}
\end{align*}
Further,
\begin{align*}
 I_1&(l,0,l+1,0)=\frac{A_l^0\cdot A_{l+1}^0}{\Sigma_n}
   \cdot\underbrace{\frac{\pi\,\Gamma(l+2\lambda+1)}{2^{2\lambda}\,l!\,(l+\lambda)\,(l+\lambda+1)\,[\Gamma(\lambda)]^2}}_{C_{\vartheta,c}(1,l,0,l+1,0)}\\
&\cdot\prod_{\iota=2}^{n-1}
   \underbrace{\frac{\pi\,\Gamma(2\lambda-\iota+1)}{2^{2\lambda-\iota}\left(\lambda+\tfrac{1-\iota}{2}\right)
      \left[\Gamma\left(\lambda+\tfrac{1-\iota}{2}\right)\right]^2}}_{C_{\vartheta,1}(\iota,0,0,0,0)}
   \cdot\underbrace{2\pi}_{C_{\varphi,1}(0)}
\end{align*}
and
\begin{align*}
 I_1&(l,2,l+1,2)=\frac{A_l^2\cdot A_{l+1}^2}{\Sigma_n}
   \cdot\underbrace{\frac{\pi\,\Gamma(l+2\lambda+3)}{2^{2\lambda+4}\,(l-2)!\,(l+\lambda)\,(l+\lambda+1)\,[\Gamma(\lambda+2)]^2}}_{C_{\vartheta,c}(1,l,2,l+1,2)}\\
&\cdot\underbrace{\frac{\pi\,\Gamma(2\lambda+1)}
      {2^{2\lambda-1}\,\left(\lambda+\tfrac{3}{2}\right)\left[\Gamma\left(\lambda-\tfrac{1}{2}\right)\right]^2}}_{C_{\vartheta,1}(2,2,0,2,0)}
   \cdot\prod_{\iota=3}^{n-1}
   \underbrace{\frac{\pi\,\Gamma(2\lambda-\iota+1)}{2^{2\lambda-\iota}\left(\lambda+\tfrac{1-\iota}{2}\right)
      \left[\Gamma\left(\lambda+\tfrac{1-\iota}{2}\right)\right]^2}}_{C_{\vartheta,1}(\iota,0,0,0,0)}
   \cdot\underbrace{2\pi}_{C_{\varphi,1}(0)}
\end{align*}
Using the doubling formula for the gamma function \cite[formula~8.335.1]{GR},
\begin{equation}\label{eq:doubling_formula}
\Gamma(2x)=\frac{2^{2x-1}}{\sqrt\pi}\,\Gamma(x)\,\Gamma\left(x+\tfrac{1}{2}\right),
\end{equation}
one can simplify the terms~$C_{\vartheta,1}$,
\begin{align*}
C_{\vartheta,1}&(\iota,0,0,0,0)
=\frac{\sqrt\pi\,\Gamma\left(\lambda-\tfrac{\iota}{2}+1\right)}{\Gamma\left(\lambda-\tfrac{\iota}{2}+\tfrac{3}{2}\right)},\\
C_{\vartheta,1}&(2,2,0,2,0)=\frac{\pi\left(\lambda-\tfrac{1}{2}\right)\Gamma(2\lambda+1)}
   {2^{2\lambda-1}\left(\lambda+\tfrac{3}{2}\right)\Gamma\left(\lambda-\tfrac{1}{2}\right)\Gamma\left(\lambda+\tfrac{1}{2}\right)}\\
&=\frac{2\sqrt\pi\left(\lambda-\tfrac{1}{2}\right)\Gamma(\lambda+1)}{\left(\lambda+\tfrac{3}{2}\right)\Gamma\left(\lambda-\tfrac{1}{2}\right)}.
\end{align*}
Consequently,
\begin{align*}
I_1&\,(l,0,l+1,0)=\frac{A_l^0\cdot A_{l+1}^0}{\Sigma_n}\cdot\frac{\pi\,\Gamma(l+2\lambda+1)}{2^{2\lambda}\,l!\,(l+\lambda)\,(l+\lambda+1)\,[\Gamma(\lambda)]^2}\\
&\cdot\pi^{\frac{n-2}{2}}\cdot\prod_{\iota=2}^{n-1}\frac{\Gamma\left(\lambda-\tfrac{\iota}{2}+1\right)}{\Gamma\left(\lambda-\tfrac{\iota}{2}+\tfrac{3}{2}\right)}\cdot2\pi\\
&=\frac{A_l^0\cdot A_{l+1}^0}{\Sigma_n}\cdot\frac{\pi^{\lambda+\frac{3}{2}}\,\Gamma(l+2\lambda+1)}{2^{2\lambda-1}\,l!\,(l+\lambda)\,(l+\lambda+1)\,[\Gamma(\lambda)]^2}
   \cdot\frac{1}{\Gamma\left(\lambda+\tfrac{1}{2}\right)}.
\end{align*}
Thus, by the doubling formula for the gamma function~\eqref{eq:doubling_formula},
\begin{align*}
I_1&(l,0,l+1,0)=\frac{A_l^0\cdot A_{l+1}^0}{\Sigma_n}
   \cdot\frac{\pi^{\lambda+1}\,\Gamma(l+2\lambda+1)}{l!\,(l+\lambda)\,(l+\lambda+1)\,\Gamma(\lambda)\,\Gamma(2\lambda)}\\
&=\frac{A_l^0\,A_{l+1}^0\,\lambda^2}{(l+\lambda)(l+\lambda+1)}\,\binom{l+2\lambda}{l}.
\end{align*}
In a similar way one obtains
\begin{align*}
 I_1(l,2,l+1,2)&=\frac{A_l^2\,A_{l+1}^2\lambda(2\lambda-1)^2(2\lambda+1)}{8(\lambda+1)(2\lambda+3)(l+\lambda)\,(l+\lambda+1)}
   \,\binom{l+2\lambda+2}{l}\cdot(l-1)l.
\end{align*}
Consequently,
\begin{align*}
\widehat G_l^0&\cdot\widehat G_{l+1}^0\cdot I_1(l,0,l+1,0)=\frac{\rho^2}{\Sigma_n^2\,\lambda^2\,A_l^0\,A_{l+1}^0}\\
&\cdot\frac{l^2(l+1)^2(l+\lambda)(l+\lambda+1)(l+2\lambda)(l+2\lambda+1)}{(2\lambda+1)^2}\cdot e^{-\rho(2l+1)}\\
&\cdot\frac{A_l^0\,A_{l+1}^0\,\lambda^2}{(l+\lambda)(l+\lambda+1)}\,\binom{l+2\lambda}{l}\\
&=\frac{\rho^2}{\Sigma_n^2(2\lambda+1)^2}\,\binom{l+2\lambda}{l}\,l^2(l+1)^2(l+2\lambda)(l+2\lambda+1)\cdot e^{-\rho(2l+1)}
\end{align*}
and
\begin{align*}
\widehat G_l^2&\cdot\widehat G_{l+1}^2\cdot I_1(l,2,l+1,2)=\frac{\rho^2}{\Sigma_n^2\,\lambda\,A_l^0\,A_{l+1}^0}
   \cdot\frac{4l^2(l+1)(l+\lambda)(l+\lambda+1)(l+2\lambda+1)}{(2\lambda+1)(4\lambda^2+8\lambda+3)}\\
&\cdot\sqrt{(l+2\lambda)(l+2\lambda+2)(l-1)(l+1)}\\
&\cdot\frac{A_l^2\,A_{l+1}^2\lambda(2\lambda-1)^2(2\lambda+1)}{8(\lambda+1)(2\lambda+3)(l+\lambda)\,(l+\lambda+1)}
   \,\binom{l+2\lambda+2}{l}\cdot(l-1)l\\
&=\frac{8\lambda(\lambda+1)\rho^2}{\Sigma_n^2(2\lambda+1)(2\lambda+3)}\cdot\binom{l+2\lambda+2}{l}(l-1)l^2(l+1)\,e^{-\rho(2l+1)}.
\end{align*}

Using notation~\eqref{eq:Smmu},
\begin{align*}
&\sum_{l=0}^\infty\widehat G_l^0\cdot\widehat G_{l+1}^0\cdot I_1(l,0,l+1,0)=\frac{\rho^2\,e^{-\rho}}{\Sigma_n^2(2\lambda+1)^2}\\
&\cdot\left[S_6^{2\lambda}+(4\lambda+3)S_5^{2\lambda}+(4\lambda^2+10\lambda+3)S_4^{2\lambda}+(8\lambda^2+8\lambda+1)S_3^{2\lambda}
   +(4\lambda^2+2\lambda)S_2^{2\lambda}\right]
\end{align*}
and
$$
\sum_{l=0}^\infty\widehat G_l^2\cdot\widehat G_{l+1}^2\cdot I_1(l,2,l+1,2)
   =\frac{8\lambda(\lambda+1)\rho^2\,e^{-\rho}}{\Sigma_n^2(2\lambda+1)(2\lambda+3)}\cdot\left[S_4^{2\lambda+2}-S_2^{2\lambda+2}\right].
$$
Note that $\widehat G_l^0\,\widehat G_{l+1}^0$ resp. $\widehat G_l^2\,\widehat G_{l+1}^2$ vanish for $l=0$ resp. $l=0$ and $l=1$ such that the range of summation in~\eqref{eq:G2_xi_O_G} can be extended to $\mathbb N_0$. Consequently, by~\eqref{eq:G2_xi_O_G} and~\eqref{eq:Smmurho}, for $\lambda\geq\frac{3}{2}$,
\begin{align*}
\|G\|^2&\cdot\xi_O(G)=\frac{\lambda+1}{\Sigma_n^2\rho^{2\lambda}}
   \cdot\left[\frac{3(\lambda+2)(\lambda+3)(2\lambda+5)}{2^{2\lambda+3}\rho^5}+\frac{3\lambda(\lambda+2)^2(2\lambda+5)}{2^{2\lambda+2}\rho^4}\right.\\
&+\frac{(\lambda+2)(24\lambda^4+80\lambda^3+48\lambda^2-31\lambda-9)}{2^{2\lambda+4}\rho^3}\\
&+\left.\frac{\lambda(8\lambda^5+32\lambda^4+24\lambda^3-31\lambda^2-40\lambda-12)}{2^{2\lambda+3}\rho^2}
   +\mathcal O\left(\frac{1}{\rho}\right)\right]\cdot\left(\begin{array}{c}1\\0\\0\\\vdots\\0\end{array}\right).
\end{align*}
This, together with~\eqref{eq:G2} yields
\begin{align*}
\xi_O(G)&=\left[1-\left(\frac{1}{2}-\frac{4}{\lambda+3}+\frac{14}{6\lambda+15}\right)\cdot\rho^2\right.\\
&+\left.\frac{2\lambda^2(3\lambda+5)}{3(\lambda+2)(\lambda+3)^2(2\lambda+5)}\cdot\rho^3+\mathcal O(\rho^4)\right]
   \cdot\left(\begin{array}{c}1\\0\\0\\\vdots\\0\end{array}\right)\qquad\text{for }\rho\to0.
\end{align*}
Hence,
\begin{align*}
\text{var}_S(G)&=\frac{6\lambda^2+13\lambda+9}{3(\lambda+3)(2\lambda+5)}\cdot\rho^2
   -\frac{4\lambda^2(3\lambda+5)}{3(\lambda+2)(\lambda+3)^2(2\lambda+5)}\cdot\rho^3\\
&+\mathcal O(\rho^4)\qquad\text{for }\rho\to0.
\end{align*}

\subsection{Localization in the momentum domain}

Using formula~\eqref{eq:varMF_series} we obtain for the variance in the momentum domain of the function~$G$ the following expression
\begin{align*}
\|G\|^2&\cdot\text{var}_M(G)=\frac{1}{\Sigma_n^2}\sum_{l=0}^\infty l(l+2\lambda)\cdot\beta_{l,0}^2(\beta_{l,0}^2+\beta_{l,1}^2)\cdot\frac{(l+\lambda)^2}{\lambda^2}
   \cdot\frac{\rho^2l^2e^{-2\rho l}}{(A_l^0)^2}\\
&=\frac{2\rho^2}{\Sigma_n^2(2\lambda+1)(2\lambda+3)}\binom{l+2\lambda}{l}l^4(l+\lambda)(l+2\lambda)(3l^2-4\lambda+6\lambda l)\,e^{-2\rho l}.
\end{align*}
It can be written as
\begin{align*}
\|G\|^2&\cdot\text{var}_M(G)=\frac{2\rho^2}{\Sigma_n^2(2\lambda+1)(2\lambda+3)}\\
&\cdot\left[3S_8^{2\lambda}+15\lambda S_7^{2\lambda}+4\lambda(6\lambda-1)S_6^{2\lambda}+12\lambda^2(\lambda-1)S_5^{2\lambda}-8\lambda^3S_4^{2\lambda}\right].
\end{align*}
Thus, by~\eqref{eq:Smmurho},
\begin{align*}
\|G\|^2&\cdot\text{var}_M(G)=\frac{(\lambda+1)(\lambda+2)(2\lambda+5)}{\Sigma_n^2\rho^{2\lambda}}
   \left[\frac{3(\lambda+3)(\lambda+4)(2\lambda+7)}{2^{2\lambda+4}\rho^7}\right.\\
&+\frac{3\lambda(\lambda+3)^2(2\lambda+7)}{2^{2\lambda+3}\rho^6}+\frac{\lambda(\lambda+3)(12\lambda^3+64\lambda^2+75\lambda-9)}{2^{2\lambda+4}\rho^5}\\
&+\frac{\lambda^2(\lambda+2)(2\lambda+3)(2\lambda^2+7\lambda-3)}{2^{2\lambda+3}\rho^4}+\left.\mathcal O\left(\frac{1}{\rho^3}\right)\right]\quad\text{for }\rho\to0
\end{align*}
if $\lambda\geq\frac{3}{2}$ and further, using asymptotics~\eqref{eq:G2},
\begin{align*}
\text{var}&_M(G)=\frac{(\lambda+4)(2\lambda+7)}{2\rho^2}+\frac{\lambda(2\lambda+7)}{(\lambda+3)\rho}\\
&-\frac{\lambda(4\lambda^4-8\lambda^3-189\lambda^2-492\lambda-351)}{6(\lambda+3)^2(2\lambda+5)}\\
&+\frac{2\lambda^2(2\lambda^4-2\lambda^3-87\lambda^2-246\lambda-189)\rho}{3(\lambda+2)(\lambda+3)^3(2\lambda+5)}+\mathcal O(\rho^2)\quad\text{for }\rho\to0.
\end{align*}

\subsection{The uncertainty product}

The square root of the product of the variances in the space and the momentum domains yields the uncertainty product of~$G$. For $\rho\to0$ it behaves like
\begin{align*}
U(G)&=\sqrt\frac{(\lambda+4)(2\lambda+7)(6\lambda^2+13\lambda+9)}{6(\lambda+3)(2\lambda+5)}\\
&-\frac{\lambda(9\lambda^2+5\lambda-18)}{(\lambda+2)(\lambda+3)}
   \sqrt\frac{2\lambda+7}{6(\lambda+3)(\lambda+4)(2\lambda+5)(6\lambda^2+13\lambda+9)}\cdot\rho
+\mathcal O(\rho^2).
\end{align*}
Note that this quantity tends to a finite value for~$\rho\to0$. In view of the discussion in~\cite[Section~4]{IIN17USW} this feature seems to be rare. Further,
$$
\lim_{\rho\to0}U(G)=\lambda+\frac{25}{12}+\mathcal O\left(\frac{1}{\lambda}\right)\qquad\text{for }\lambda\to\infty,
$$
i.e., if the dimension of the sphere tends to infinity, $\lim_{\rho\to0}U(G)$ approaches the value $\lambda+\frac{25}{12}$, close to the optimal one $\lambda+\frac{1}{2}$, compare formula~\eqref{eq:UP}.

\begin{figure}[t]\label{fig:proportion}
\includegraphics[width=0.8\linewidth]{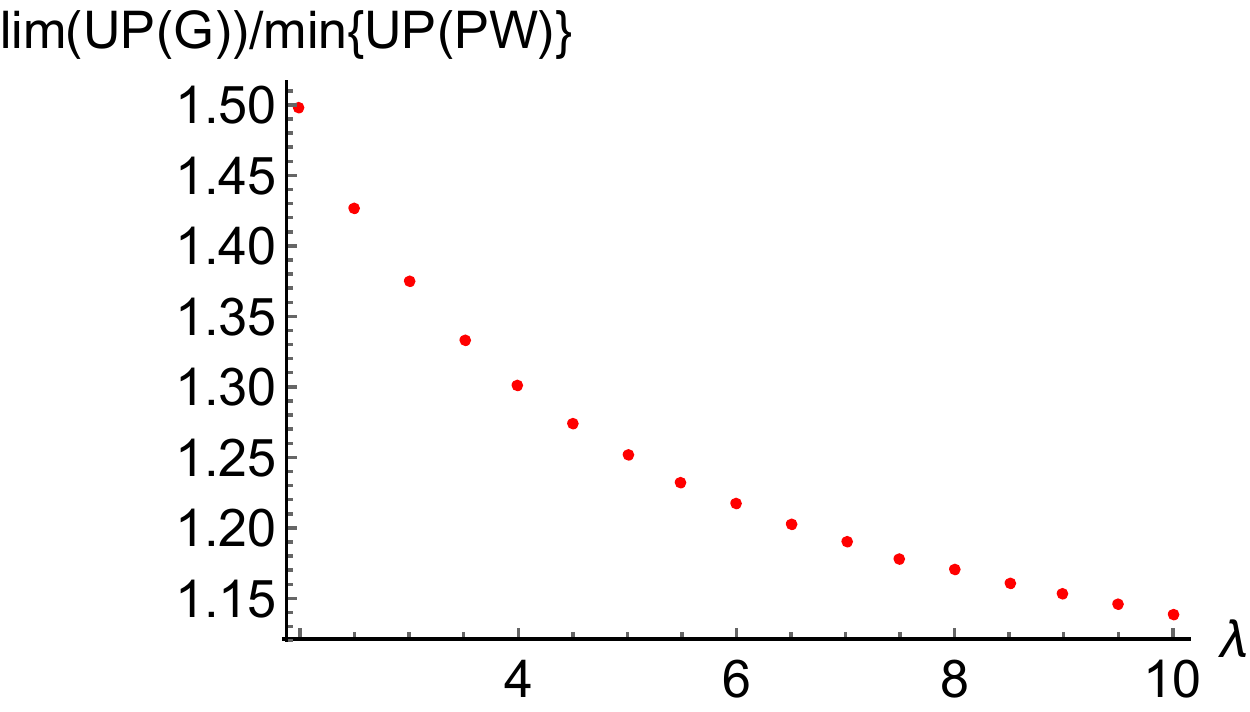}
\caption{Proportion of~$\lim_{\rho\to0}U(G)$ to the minimum limiting value of the uncertainty product of the zonal Poisson wavelets}
\end{figure}

Figure~\ref{fig:proportion} illustrates the ratio of the limiting value of the uncertainty product of~$G$ and the minimum of the limiting values of the uncertainty product of the zonal Poisson wavelets~$g_\rho^m$, which is obtained for $m=[\lambda]$, see the discussion in~\cite[Section~3]{IIN17UPW}. The considered dimensions of the sphere are greater than or equal to~$5$ because of the exceptional values of the uncertainty product of the Poisson wavelets in dimensions~$3$ and~$4$, compare \cite[Theorem~3.1]{IIN17UPW}. The discussed ratio is monotonously decreasing for $\lambda\geq2$ (all of the roots of the derivative of its square are less than~$1$) and tends to~$1$ for~$\lambda$ tending to infinity.


\begin{thebibliography}{99}
\bibitem{sB09} S.~Bernstein, \emph{Spherical singular integrals, monogenic kernels and wavelets on the three--dimensional sphere}, Adv. Appl. Clifford Algebr. 19 (2009), no.~2, 173--189.
\bibitem{BE10KBW} S. Bernstein and S. Ebert, \emph{Kernel based wavelets on $S^3$}, J. Concr. Appl. Math. 8 (2010), no. 1, 110--124.
\bibitem{BE10WS3} S. Bernstein and S. Ebert, \emph{Wavelets on $S^3$ and $SO(3)$ --- their construction, relation to each other and Radon transform of wavelets on $SO(3)$}, Math. Methods Appl. Sci. 33 (2010), no. 16, 1895--1909.
\bibitem{CPMDHJ05} A. Chambodut, I. Panet, M. Mandea, M. Diament, M. Holschneider, and O.~Jamet, \emph{Wavelet frames: an alternative to spherical harmonic representation of potential fields}, Geophys. J. Int. 163 (2005), 875--899.
\bibitem{DX15b} F. Dai and Y. Xu, \emph{Erratum to: The Hardy-Rellich inequality and uncertainty principle on the sphere}, Constr. Approx. 42 (2015), no. 1, 181--182. 
\bibitem{DX15a} F. Dai and Y. Xu, \emph{The Hardy-Rellich inequality and uncertainty principle on the sphere}, Constr. Approx. 40 (2014), no. 1, 141--171.
\bibitem{DQC17} P. Dang and T. Qian, and Q. Chen, \emph{Uncertainty principle and phase-amplitude analysis of signals on the unit sphere}, Adv. Appl. Clifford Algebr. 27 (2017), no. 4, 2985--3013.
\bibitem{EBCK09} S. Ebert, S. Bernstein, P. Cerejeiras, and U. K\'ahler, \emph{Nonzonal wavelets on~$\mathcal S^N$}, 18$^\text{th}$ International Conference on the Application of Computer Science and Mathematics in Architecture and Civil Engineering, Weimar 2009.
\bibitem{FGS-book} W. Freeden, T. Gervens, and M. Schreiner, \emph{Constructive approximation on the sphere. With applications to geomathematics}, Numerical Mathematics and Scientific Computation, The Clarendon Press, Oxford University Press, New York, 1998.
\bibitem{FW-C97} W.~Freeden and U.~Windheuser, \emph{Combined spherical harmonic and wavelet expansion -- a future concept in {E}arth's gravitational determination}, Appl. Comput. Harmon. Anal. 4 (1997), no. 1, 1--37.
\bibitem{FW96} W.~Freeden and U.~Windheuser, \emph{Spherical wavelet transform and its discretization}, Adv. Comput. Math. 5 (1996), no.~1, 51--94.
\bibitem{GG04a}  S.S. Goh and T.N.T. Goodman, \emph{Uncertainty principles and asymptotic behavior}, Appl. Comput. Harmon. Anal. 16 (2004), no. 1, 19--43.
\bibitem{GG04b} T.N.T. Goodman and S.S. Goh, \emph{Uncertainty principles and optimality on circles and spheres}, Advances in constructive approximation: Vanderbilt 2003, 207--218,
Mod. Methods Math., Nashboro Press, Brentwood, TN, 2004. 
\bibitem{GR} I.S. Gradshteyn and I.M.~Ryzhik, \emph{Table of integrals, series, and products}, Elsevier/Academic Press, Amsterdam, 2007.
\bibitem{HH09} M.~Hayn and M.~Holschneider, \emph{Directional spherical multipole wavelets}, J. Math. Phys. 50 (2009), no. 7, 073512, 11 pp.
\bibitem{mH96} M.~Holschneider, \emph{Continuous wavelet transforms on the sphere}, J. Math. Phys. 37 (1996), no. 8, 4156--4165.
\bibitem{HI07} M.~Holschneider and I.~Iglewska-Nowak, \emph{Poisson wavelets on the sphere}, J. Fourier Anal. Appl. 13 (2007), no. 4, 405--419.
\bibitem{IIN15CWT} I.~Iglewska-Nowak, \emph{Continuous wavelet transforms on n-dimensional spheres}, Appl. Comput. Harmon. Anal. 39 (2015), no. 2, 248--276.
\bibitem{IIN18DW} I.~Iglewska-Nowak, \emph{Directional wavelets on n-dimensional spheres}, Appl. Comput. Harmon. Anal. 44 (2018), no. 2, 201--229.
\bibitem{IIN17FDW} I.~Iglewska-Nowak, \emph{Frames of directional wavelets on n-dimensional spheres}, Appl. Comput. Harmon. Anal. 43 (2017), no. 1, 148--161.
\bibitem{IIN16MR} I.~Iglewska-Nowak, \emph{Multiresolution on n-dimensional spheres}, Kyushu J. Math. 70 (2016), no. 2, 353--374.
\bibitem{IIN17USW} I.~Iglewska-Nowak, \emph{On the uncertainty product of spherical wavelets}, Kyushu J. Math. 71 (2017), no. 2, 407--416.
\bibitem{IIN15PW} I.~Iglewska-Nowak, \emph{Poisson wavelets on n-dimensional spheres}, J. Fourier Anal. Appl. 21 (2015), no. 1, 206--227.
\bibitem{IIN17UPW} I.~Iglewska-Nowak, \emph{Uncertainty of Poisson wavelets}, Kyushu J. Math. 71 (2017), 349--362.
\bibitem{nLF07} N. La\'in Fernández, \emph{Optimally space-localized band-limited wavelets on $\mathfrak S^{q-1}$}, J. Comput. Appl. Math. 199 (2007), no. 1, 68--79.
\bibitem{nLF03} N. La\'in Fern\'andez, \emph{Polynomial bases on the sphere}, doctoral thesis, Universit\"at zu L\"ubeck 2003.
\bibitem{LFP02} N. Laín Fernández, N. and J. Prestin, \emph{Localization of the spherical Gauss-Weierstrass kernel}, Constructive theory of functions, 267--274, DARBA, Sofia, 2003. 
\bibitem{NW96} F. J. Narcowich and J. D. Ward, \emph{Nonstationary wavelets on the m-sphere for scattered data}. Appl. Comput. Harmon. Anal. 3 (1996), no. 4, 324--336.
\bibitem{RV97} M. R\"osler and M. Voit, \emph{An uncertainty principle for ultraspherical expansions}, J. Math. Anal. Appl. 209 (1997), no. 2, 624--634.
\bibitem{kS02} K.~Selig, \emph{Uncertainty principles revisited}, Orthogonal polynomials, approximation theory, and harmonic analysis (Inzel, 2000). Electron. Trans. Numer. Anal. 14 (2002), 165--177.
\bibitem{sS15} S. Steinerberger, \emph{An uncertainty principle on compact manifolds}, J. Fourier Anal. Appl. 21 (2015), no. 3, 575--599.
\bibitem{Vilenkin} N. Ja. Vilenkin, \emph{Special functions and the theory of group representations}, in Translations of Mathematical Monographs, Vol. 22, American Mathematical Society, Providence, R. I., 1968.

\end{thebibliography}
\end{document}